\documentclass[final,5p,times,twocolumn]{elsarticle}

\pdfoutput=1
\usepackage{amssymb}
\usepackage{amsfonts}
\usepackage{latexsym, amssymb,amsmath, amsbsy,amsopn, amstext}
\usepackage{amsmath, amsbsy}
\usepackage{hyperref}
\usepackage{array,hhline}
\usepackage{makecell}
\usepackage{float}
\usepackage{ifpdf}
\usepackage{xcolor}
\usepackage{tikz}
\usepackage{algorithm}
\usetikzlibrary{arrows,shapes,snakes,shadows,positioning,automata,patterns}
\usetikzlibrary{trees,decorations.pathmorphing,decorations.markings}
\usepackage{bm}

\def\argmin{\mathop{\rm argmin}}
\newtheorem{theorem}{Theorem}
\newtheorem{lemma}{Lemma}

\newtheorem{assumption}{Assumption}
\newtheorem{proposition}{Proposition}
\newtheorem{cor}{Corollary}
\newtheorem{remark}{Remark}

\newcommand{\blue}[1]{{\color{black}#1}}

\allowdisplaybreaks

\journal{}

\begin{document}

\begin{frontmatter}



\title{An Optimal Distributed Algorithm with Operator Extrapolation for Stochastic Aggregative Games}


%
%
\tnotetext[t1]{This research is supported by  Shanghai Municipal Science and Technology Major Project No. 2021SHZDZX0100,  Shanghai Municipal Commission of Science and Technology Project No. 19511132101,
National Natural Science Foundation of China under Grant 62003239, Shanghai Sailing Program No. 20YF1453000,
 and the Fundamental Research Funds for the Central Universities, China.}
\author[rvt,gf]{Tongyu Wang}
\ead{wang\_tongyu@tongji.edu.cn}
\author[rvt,ai,focal]{Peng Yi }
\ead{yipeng@tongji.edu.cn}
\author[rvt,ai,fsc]{Jie Chen }
\ead{chenjie206@tongji.edu.cn}
\ead[url]{http://www.elsevier.com}

\address[rvt]{Department of Control Science and Engineering, Tongji University,  Shanghai,   201804, China;}
\address[ai]{The Shanghai Research
Institute for Intelligent Autonomous Systems,Shanghai, 201804, China.}
\address[focal]{The Institute of Advanced Study, Tongji University,  Shanghai, 201804,  China}
\address[fsc]{Frontiers Science Center for Intelligent Autonomous Systems, Tongji University, Shanghai, China}
\address[gf]{Shanghai Institute of Intelligent Science and Technology, Tongji University, Shanghai 200092, China}


\begin{abstract}
This work studies  Nash equilibrium seeking for a class of stochastic aggregative games,  where each player has an expectation-valued objective function depending on its local strategy and the aggregate of all players' strategies. We propose a distributed algorithm with operator extrapolation, in which each player maintains an estimate of this aggregate by exchanging this information with its neighbors over a time-varying network, and updates its decision through the mirror descent method. An operator extrapolation at the search direction is applied such that the two step historical gradient samples are utilized  to accelerate the convergence.  Under the strongly monotone assumption on the pseudo-gradient mapping, we prove that  the proposed algorithm can achieve the optimal convergence rate of $\mathcal{O}(1/k)$ for Nash equilibrium seeking of stochastic games.  Finally, the algorithm  performance is demonstrated via numerical simulations.
\end{abstract}


\end{frontmatter}


\section{Introduction}

In recent years, the non-cooperative games have been widely applied in decision-making for networked systems, where selfish networked players aim to optimize their own objective functions that have coupling with other players' decisions.
Nash equilibrium (NE) is one of the  most widely used solution concept  to non-cooperative games, under which no participant can improve its  utility by unilaterally deviating from the equilibrium strategy \cite{nash50equilibrium}. 
Specifically,   the  aggregative games is an important class  of non-cooperative games, in which each player's cost function   depends on its own strategy and the aggregation of   all player strategies \cite{jensen2010aggregative}.
This aggregating feature emerges in numerous decision-making problems over networks systems, when the individual utility is affected by network-wide aggregate.  Hence,  aggregative games and its NE seeking algorithm get wide research, with application in charging coordination for plug-in electric vehicles \cite{grammatico2015decentralized,grammatico2017dynamic}, multidimensional opinion dynamics\cite{GHADERI20143209,parise2020distributed}, communication channel congestion control \cite{Bakhshayesh2021}, and energy management in power grids \cite{belgioioso2021operationallysafe,DengAggregate2021}, etc. It also worth pointing out that stochastic game should be taken as a proper decision-making model when there are information uncertainties, dating back to \cite{shapley1953stochastic}, while the NE seeking for stochastic games with various sampling schemes are recently studied \cite{lei2020synchronous,lei2018linearly}. Hence, when the game has both information uncertainties and  aggregating-type objective functions, stochastic
aggregative game and its NE seeking are studied \cite{meigs2019learning,lei2021distributed,franci2022stochastic}.

Generally speaking, NE seeking methods include centralized, distributed, and semi-decentralized methods. In the centralized methods \cite{Frihauf2012}, each player  needs to know the complete  information  to compute an NE in an introspective manner.
In the semi-decentralized  methods \cite{belgioioso2017semi},   there exists a central coordinator  gathering and broadcasting signals   to all   players. While  in   distributed methods \cite{Ye2017}, each player has its private information   and   seeks the NE through local communications with   other players.
As such, designing distributed NE seeking algorithm has gained tremendous research interests due to the benefits owning to decentralization.

Existing works on distributed NE seeking largely resided in best-response(BR) schemes \cite{shanbhag16inexact}\cite{lei2018linearly}\cite{parise2020distributed} and gradient-based approaches \cite{salehisadaghiani2018distributed}\cite{yi2019operator}\cite{Zhou2017}. In BR schemes, players select a strategy that best responds to  its rivals' current strategies. For example, \cite{lei2019synchronous} proposed an inexact generalization of this scheme in which an inexact best response strategy is computed via a stochastic approximation scheme. In gradient-based approaches, the algorithm is  easily implementable  with a low computation complexity per iteration. For instance, \cite{Tatarenko2021} designed an accelerated distributed gradient play for a strongly monotone game  and proved its geometric  convergence, while it needs to estimate  all players' strategies.
To reduce the communication costs, \cite{koshal2016distributed} proposed a consensus-based distributed gradient algorithm for aggregative games, in which players only need to estimate the aggregate in each iteration. Furthermore, \cite{lei2020Agg} proposed a distributed gradient algorithm based on iterative Tikhonov regularization method to resolve the class of monotone aggregative games. \cite{ANAG2022} proposed a fully asynchronous distributed algorithm and rigorously show the convergence to a Nash equilibrium.
Besides, \cite{liang2017distributed,DengAggregate2021} proposed a continuous-time distributed algorithm for aggregative games.

\begin{table*}[htp]
	\centering
	\setlength{\tabcolsep}{5mm}{
		\begin{tabular}{|l|l|l|l|l|l|l|}
			\hline
			\makecell[c]{Literature} & \makecell[c]{stochastic}
			& \makecell[c]{Method}  & \makecell[c]{rate}
			& \makecell[c]{Iteration} & \makecell[c]{communication} &sample   \\
			\hline
		    \makecell[c]{\cite{parise2020distributed}} &\makecell[c]{$\times$} & \makecell[c]{BR-based} & \makecell[c]{-}
		    &  \makecell[c]{-} & \makecell[c]{$v$ rounds}  & \makecell[c]{-} \\
			\hline
		    \makecell[c]{\cite{ANAG2022}} &\makecell[c]{$\times$} & \makecell[c]{gradient-based}  & \makecell[c]{$1/k$}
		    &  \makecell[c]{-} 	& \makecell[c]{$1$ round} & \makecell[c]{-}\\
			\hline
			\makecell[c]{\cite{Lei2018agg}}   & \makecell[c]{\checkmark}  & gradient-based  & \makecell[c]{$1/k$}
			& \makecell[c]{$\mathcal{O}(\ln (1 / \epsilon))$}
			& \makecell[c]{$k+1$ rounds}    & \makecell[c]{$N_{k}$}\\
			\hline
			This work & \makecell[c]{\checkmark}
			&\makecell[c]{gradient-based}  & \makecell[c]{$1/k$}
			& \makecell[c]{Corollary\ref{cor1}}
			& \makecell[c]{$1$ round} & \makecell[c]{1}     \\
			\hline
	\end{tabular}}
\caption{A list of some recent research papers on distributed  methods for resolving aggregative games.}\label{tab1}
In column 2 ``stochastic", $\checkmark$ means that the studied problem  is stochastic, while $\times$ implies that the studied problem is deterministic; In columns 4 and 5, dash - implies that the literature has not studied this property; In column 6 ``communication" stands for that the communication round   required at each iteration; In column 7 ``sample" stands for the  number of sampled gradients per iteration.
\end{table*}

Moreover, fast convergence rate is indispensable bonus for distributed Nash equilibrium seeking algorithm, since it implies
less communication rounds. Particulary, for stochastic  games,  distributed NE seeking with a fast convergence rate is highly desirable since its also implies low sampling cost of stochastic gradients.
Motivated by this, \cite{enrich2019extragradient} designed an extra-gradient method to improve the speed of NE seeking, while it needs two steps of operator evaluations and two step of projections at each iteration. To further reduce the computation burden, \cite{Malitsky2015} proposed  a simpler recursion with one step of operator evaluation and  projection at each iteration, given by $x_{t+1}=\operatorname{argmin}_{x \in X} \gamma_{t}\left\langle F\left(x_{t}+\beta_{t}\left(x_{t}-x_{t-1}\right)\right),x\rangle +V\left(x_{t}, x\right)\right.$, called projected reflected gradient method. But the iterate $x_{t}+\beta_{t}\left(x_{t}-x_{t-1}\right)$ may sit outside the feasible set $X$,
 so it needs to impose the strong monotone assumption on the pseudo-gradient over $\mathbb{R}^{n}$. Recently, \cite{kotsalis2020simple} proposed an operator extrapolation (OE) method to solve the stochastic variational inequality problem,  for which  the optimal convergence rate can be achieved through  one operator evaluation and a single  projection per iteration.

Motivated by the demands for fast NE seeking algorithm of stochastic aggregative games and inspired by \cite{kotsalis2020simple}, we  propose a distributed operator extrapolation (OE) algorithm via mirror descent and dynamical averaging tracking. At each stage, each player aligns its intermediate estimate by a consensus step with  its neighbors' estimates of the aggregate, samples an unbiased estimate of its payoff gradients, takes a small step via the OE method, and then mirrors it back to the strategy set.
  The algorithm achieves the optimal convergence rate $\mathcal{O}(1/k)$ for the class of stochastic strongly monotone aggregative games. The numerical experiments on a Nash-cournot competition problem demonstrate it advantages over some existing NE seeking methods.
 In addition, we compare our algorithm with some other  distributed algorithms for  aggregative games and  summarize those works in Table \ref{tab1}.
  Specially, we compare the iteration complexity and communication rounds to achieve an $\epsilon-\mathrm{NE}$  when the iterate $x$ satisfies $\mathbb{E}[\| x-$ $\left.x^{*} \|^{2}\right] \leq \epsilon$.
The proposed algorithm can achieve the optimal convergence rate with only one communication round   per  iteration in comparison with multiple communication rounds required by \cite{parise2020distributed} and \cite{Lei2018agg}. Though \cite{ANAG2022} can achieve  convergence  rate $1/k$ with fixed step-sizes for a deterministic game,   our work can achieve the optimal convergence result with diminishing step for a stochastic game, which can model various information uncertainties.

The  rest of paper is organized as follows. In section \ref{sec:formulation}, we give the formulation of  the stochastic aggregate games, and state the assumptions. In section \ref{sec:agg}, we propose a distributed operator extrapolation algorithm    and provide convergence results for the class of strongly monotone games.
The proof of main results are given in Section \ref{sec:proof}. We  show the empirical performance of our algorithm through the Nash-cournot  model in section \ref{sec:sim},
and give concluding remarks in Section \ref{sec:con}.

{\bf Notations:} A vector $x$ is a column vector while $x^T$ denotes its transpose. $\langle  x,y\rangle=x^Ty$ denotes the inner product of vectors $x,y.$  $\|x\|$
{denotes} the Euclidean vector norm, i.e., $\|x\|=\sqrt{x^Tx}$. A nonnegative square   matrix $A $   is  called doubly stochastic if
$A\mathbf{1} =\mathbf{1}$ and  $\mathbf{1}^T A =\mathbf{1}^T$, where
$\mathbf{1}$ denotes the vector with each entry equal 1.
$\mathbf{I}_N \in \mathbb{R}^{N\times N}$ denotes the  identity matrix.
Let $\mathcal{G}=\{ \mathcal{N}, \mathcal{E}\}$ be  a directed graph   with
$\mathcal{N}=\{1,\cdots,N\} $  denoting the set of  players and  $\mathcal{E}$
denoting  the set of  directed edges between players, where
$(j,i)\in\mathcal{E }$ if  player  $i$ can  obtain  information from
player  $j$.
The graph  $\mathcal{G}$ is  called  strongly connected if  for
any $  i,j\in \mathcal{N}$ there exists a directed path from  $i$ to
$j$, i.e., there exists  a  sequence  of edges   $ (i,i_1),(i_1,i_2),\cdots,(i_{p-1},j)$ in the digraph  with  distinct  nodes $ i_m \in \mathcal{N},~\forall m: 1 \leq m \leq p-1$.

\section{Problem  Formulation }\label{sec:formulation}

In this section, we formulate the stochastic aggregative games over networks and introduce basic assumptions.

\subsection{Problem Statement}

We consider a set of $N$ non-cooperative players indexed by $\mathcal{N} \triangleq \{1, \cdots,N\} $. Each player $i\in{\mathcal{N}}$ choose its strategy $x_{i}$ from a strategy set $X_{i} \in \mathbb{R}^{m}$. Denote by   $\mathbf{x}\triangleq \blue{(x_1^T,\cdots,x_N^T)^T} \in
\mathbb{R}^{  mN}$ and     $\mathbf{x}_{-i}\triangleq \{x_j\}_{j
	\neq i}   $   the strategy profile and the rival
strategies,  respectively. In an aggregative game, each player $i$ aim to minimize its cost function $f_i(x_i,\sigma(x))$, where  $ \sigma(\mathbf{x})   \triangleq   \sum_{j=1}^N  x_j $ is an aggregate of all players' strategies. Furthermore, given $\sigma(\mathbf{x}_{-i})\triangleq  \sum_{j=1,j\neq i}^N x_j$, the objective of player $i$ is to minimize its parameterized stochastic local cost:
\begin{align} \label{Ngame_agg} \min_{x_i \in X_i} f_i(x_i,x_i+\sigma(\mathbf{x}_{-i}))
\triangleq \mathbb{E}\left[\psi_i(x_i,x_i+\sigma(\mathbf{x}_{-i}) ; {\xi_i(\omega)}) \right], \end{align}
where  $ {\xi_i}: { \Omega} \to  {\mathbb{R}^{d_i}}$ is  a random variable  defined on the probability space $( \Omega, \mathcal{F}, \mathbb{P}_i)$,    $\psi_i : \mathbb{R}^{m} \times   \mathbb{R}^{m} \times   \mathbb{R}^{d_i}
\to \mathbb{R}$ is a scalar-valued  function, and $\mathbb{E}[\cdot]$
denotes the expectation with respect to the probability measure $ \mathbb{P}_i $.

A Nash equilibrium (NE) of  the  aggregative game
\eqref{Ngame_agg} is a  tuple $\mathbf{x}^*=\{x_i^*\}_{i=1}^N  $  such that  for each $ i  \in \mathcal{ N}$,
\begin{align*} f_i(x_i^*, \sigma(\mathbf{x}^*)) \leq
f_i(x_i, x_i +\sigma(\mathbf{x}_{-i}^*) ),  \quad \forall x_i\in X_i.  \end{align*}

We consider the scenario  that each player $i\in \mathcal{N}$  knows the structure of its private function $f_i$,  but have no  access to the aggregate   $\sigma(\mathbf{x}).$ Instead, each player may   communicate with its neighbors
over a time varying graph $\mathcal{G }_k = \{ \mathcal{N }, \mathcal{E }_k\}$.
Define $W_k  =[ \omega_{ij,k}]_{i,j =1}^N$ as the adjacency matrix, where  $ \omega_{ij,k}>0$  if  and only if  $(j,i )\in \mathcal{E }_k$,   and  $ \omega_{ij,k}=0$, otherwise. Denote  by $ N_{i,k} \triangleq \{ j \in \mathcal{N}:(j,i )\in \mathcal{E}_k \}$   the   neighboring set  of player  $i$  at time $k$.

\subsection{Assumptions}

We impose the following conditions on  the time-varying  communication graphs   $\mathcal{G }_k = \{ \mathcal{N }, \mathcal{E }_k\}$.
\begin{assumption}~\label{ass-graph} (a)   $W_k $ is   doubly stochastic for any $k\geq 0$;
	\\(b) There exists a constant  $0< \eta<1$  such that
	$$ \omega_{ij,k} \geq \eta  , \quad  \forall j \in  \mathcal{N}_{i,k},~~\forall i \in \mathcal{N},~\forall k \geq 0;$$
	(c)   There exists  a positive integer $B $  such that the union graph
	$\big \{ \mathcal{N }, \bigcup_{l=1}^B \mathcal{E }_{k+l} \big \}$ is strongly connected  for all $k\geq 0$.
\end{assumption}

We define a transition matrix $\Phi(k,s)  =W_kW_{k-1}\cdots W_s$ for any $k\geq s\geq 0$ with $\Phi(k,k+1)  =\mathbf{I}_N$,
and state a result that will be used in the sequel.

\begin{lemma}\cite[Proposition 1]{nedic2009distributed}
	Let  Assumption \ref{ass-graph} hold. Then there exist  $\theta=(1-\eta/(4N^2))^{-2}>0$ and $\beta =(1-\eta/(4N^2))^{1/B}  $ such that for any $k\geq s \geq 0,$
	\begin{align}\label{geometric}
	\left | \left[\Phi(k,s)\right]_{ij}-1/N \right| \leq \theta \beta^{k-s},\quad \forall i,j\in \mathcal{N}.
	\end{align}
\end{lemma}

We require the player-specific problem to be convex and  continuously differentiable.

\begin{assumption}~\label{ass-payoff}  For each player $i\in \mathcal{N},$
	(a) the strategy set  $X_i$ is   closed, convex and compact;
	(b) the cost function $ f_i(x_i,\sigma(x ))$ is convex  in $x_i\in X_i$
	for every fixed  $x_{-i} \in X_{-i} \triangleq \prod_{j \neq i} X_j$;
	(c)  $ f_i(x_i,\sigma )$ is continuously differentiable
	in $(x_i,\sigma) \in X_i \times  \mathbb{R}^m$.
\end{assumption}

For any $x\in X\triangleq \prod_{i=1}^N X_i$  and   $z\in \mathbb{R}^m$, define
\begin{align}\label{def-phi}
&   F_i(x_i, z) \triangleq \big(\nabla_{x_i}f_i(\cdot,\sigma) +  \nabla_{\sigma}f_i(x_i,\cdot)\big)\mid_{\sigma = z}, \notag
\\&\phi_i(\mathbf{x}) \triangleq \nabla_{x_i}f_i(x_i,\sigma(\mathbf{x}))=F_i(x_i,  \sigma(\mathbf{x})),\mbox{~and~}
\phi(\mathbf{x})\triangleq  \left(\phi_i(\mathbf{x}) \right)_{i=1}^N.
\end{align}

From \eqref{def-phi}  and   Assumption \ref{ass-payoff}(c),
it follows that  the  pseudogradient $\phi(x)$  is continuous.
Since each player-specific  problem is convex,  by \cite[Proposition 1.4.2]{facchinei02finite}, $\mathbf{x}^*$ is a NE  of \eqref{Ngame_agg} if and only if $\mathbf{x}^*$ is a solution  to a variational inequality problem VI$(X, \phi)$, i.e.,  finding
$\mathbf{x}^* \in X$ such that  \begin{align} \label{VI}
(\tilde{\mathbf{x}} - \mathbf{x}^*)^{\top} \phi(\mathbf{x}^*) \, \geq \, 0, \quad \forall  \tilde{\mathbf{x}} \in X\triangleq \prod_{j=1}^N X_j.
\end{align}
Since $\phi$ is  continuous, $X$ is   convex and compact, the  existence of NE follows immediately by \cite[Corollary 2.2.5]{facchinei02finite}.

For each $i\in \mathcal{N},$ define
\begin{align}
& M_i \triangleq   \max_{x_i \in X_i} \|x_i\|,~ M_H\triangleq \sum_{j=1}^N M_i, \label{def-bdst} \\
 &\tilde{C}\triangleq  \theta   M_H  +  2 \theta    \sum_{j=1}^N  M_j/(1-\beta). \label{def-tC}
\end{align}
We impose the following Lipschitz continuous conditions.
\begin{assumption}\label{ass-hF}
	(a)  $ \phi (\mathbf{x} )  $  is   $L$-Lipschitz continuous in $x $, i.e.,
	\[ \|\phi(\mathbf{x} ) -\phi(\mathbf{y})\| \leq L\|x-y\|,\quad \forall \mathbf{x},\mathbf{y} \in X. \]
	(b)  For each   $i\in \mathcal{N} $ and any   $x_i\in X_i$,  $ F_i(x_i, z)  $  is    $L_{fi}$-Lipschitz continuous in $z $ over some compact set, particularly,
	for all $  z_1,z_2\in  {\mathbb{R}^d}$ with $\|z_1\|\leq N\tilde{C}+M_H$ and $\|z_2\|\leq N\tilde{C}+M_H$: $$\|  F_i(x_i, z_1)- F_i(x_i, z_2)\| \leq L_{fi} \|z_1-z_2\|  .$$
\end{assumption}

%

In addition, we require  each  $  \psi_i(x_i, \sigma ;\xi_i)$  to be  differentiable, and  assume there exists a stochastic oracle  that returns  unbiased gradient sample with bounded variance.

\begin{assumption}\label{ass-noise}For each player $i\in \mathcal{N} $ and any   $ {\xi_i\in\mathbb{R}^{d_i}}$, \\
	(a)  $  \psi_i(x_i, \sigma;\xi_i)$ is  differentiable  in $x_i \in X_i$  and $\sigma \in \mathbb{R}^m $;\\
	(b)   for any  $x_i\in X_i $ and $ z \in \mathbb{R}^m $,
	$ q_i(x_i,z ;\xi_i) = (\nabla_{x_i}\psi_i(x_i, \sigma;\xi_i) + \nabla_{\sigma}  \psi_i(x_i, \sigma;\xi_i))\mid_{\sigma = z }$ satisfies
	$ \mathbb{E}\left[ q_i(x_i,z_i;\xi_i) | x_i \right] =F_i(x_i,z_i)  $
	and $  \mathbb{E}[\|  q_i(x_i,z_i ;\xi_i) -F_i(x_i,z_i)\|^2 |x_i]\leq  \nu_i^2 $ for some constant $\nu_i>0.$
\end{assumption}

\section{Algorithm Design and Main Results }\label{sec:agg}

In this section, we   design a    distributed   NE seeking algorithm  and
prove its optimal convergence  to the    Nash equilibrium in the mean-squared sense.

\subsection{Distributed  Mirror Descent Algorithm with Operator Extrapolation}


We assume throughout that the paper that the regularization function
$h$ is $1$-strongly convex, i.e.,
\[h( y) \geq   h(x ) +\langle \nabla h(x), y-x \rangle+ {1\over 2} \|x-y\|^2,
\quad \forall x,y\in \mathbb{R}^m.\]
We define the    Bregman divergence  associated with   $h$ as follows.
\begin{equation}\label{def_breg}
D (x ,y) \triangleq   h(y ) - (h(x )+ \langle\nabla h(x),y-x) \rangle, \quad \forall x,y\in \mathbb{R}^m.
\end{equation}

Recall that the operator extrapolation method for VI in \cite{kotsalis2020simple} requires a simple recursion at  each iteration given by
\[ x_{t+1}=\operatorname{argmin}_{x \in X} \gamma_{t}\left\langle F\left(x_{t}\right)+\lambda_{t}\left(F\left(x_{t}\right)-F\left(x_{t-1}\right)\right), x\right\rangle+D\left(x_{t}, x\right),\] which only involves one operator evaluation $F\left(x_{t}\right)$ and one prox-mapping over the set $X$.

Suppose that each player $i$ at   stage $k=1,2,\cdots $ selects a strategy  $x_{i,k} \in X_i $ as an estimate of its equilibrium strategy,
and  holds an estimate $v_{i,k}$ for the average  aggregate.
At  stage  $k+1$, player $i$   observes its neighbors' past information $ v_{j,k} ,j\in \mathcal{N}_{i,k}$ and updates an
intermediate estimate by the consensus step \eqref{alg-consensus},
then it computes its partial  gradient based on the sample observation,   and updates  its strategy $x_{i,k+1} $ by a mirror descent scheme  \eqref{alg-strategy0} with operator extrapolation \cite{kotsalis2020simple} by setting  $q_i(x_{i,0}, N\hat{v}_{i,1};\xi_{i, 0})=q_i(x_{i,1}, N\hat{v}_{i,2};\xi_{i,1})  $ without loss of generality.
Finally, player $i$ updates the average aggregate with the renewed strategy  $x_{i,k+1} $ by
the dynamic average  tracking scheme  \eqref{alg-average}. The procedures are summarized in  Algorithm \ref{alg1}.

\begin{algorithm}[htbp]
	\caption{Distributed    Mirror Descent Method with Operator Extrapolation}  \label{alg1}
	{\it Initialize:} Set $k=1,$ $x_{ i,1}  \in X_i$  and $ v_{ i,1}= x_{i,1}  $ for  each $i \in\mathcal{N}$.
	
	{\it Iterate until convergence}\\
	{\bf Consensus.} Each  player computes an intermediate estimate by
	\begin{align}\label{alg-consensus}
	\hat{v}_{i,k+1} = \sum_{j\in \mathcal{N}_{i,k}} w_{ij,k} v_{j,k}  .
	\end{align}
	{\bf Strategy Update.} Each  player $ i\in \mathcal{N}$ updates its    equilibrium strategy   and its estimate of the average aggregate by
	\begin{align}
	 x_{i,k+1}  &=\argmin_{x_i \in X_i} \Big( \alpha_k \big \langle
	 (1+   \lambda_k ) q_i(x_{i,k}, N\hat{v}_{i,k+1};\xi_{i,k}) \notag
 \\& - \lambda_k q_i(x_{i,k-1}, N\hat{v}_{i,k };\xi_{i,k-1})  , x_i\big \rangle + D ( x_{i,k},x_i) \Big),\label{alg-strategy0}
   \\ v_{i,k+1} &= \hat{v}_{i,k+1}+x_{i,k+1}-x_{i,k}\label{alg-average}
	\end{align}
	where $\alpha_k>0,  \lambda_k >0$, and $\xi_{i,k}$ denotes a   random  realization of $\xi_i$ at time $k$.
\end{algorithm}

Define the gradient noise
\begin{align}\label{def-zeta}
\zeta_{i,k}\triangleq q_i(x_{i,k}, N\hat{v}_{i,k+1};\xi_{i,k})-F_i(x_{i,k}, N\hat{v}_{i,k+1} ),
\end{align}
$ \mathbf{x}_k=(x_{1,k}^T,\cdots,x_{N,k}^T)^T,$
and $\mathcal{F}_k\triangleq \{\mathbf{x}_1, \xi_{i,l}, i\in \mathcal{N},  l=1,2,\cdots, k-1\}.$  Then with Algorithm \ref{alg1}, $x_{i,k }$ and $\hat{v}_{i,k+1}  $ are adapted to $\mathcal{F}_k.$ Define $\bm{\zeta}_k=(\zeta_{1,k}^T,\cdots,\zeta_{N,k}^T)^T.$
From  Assumption \ref{ass-noise}, it follows  that for each $i\in \mathcal{N}:$
\begin{equation}\label{bd-noise}
\begin{split}
&\mathbb{E}[\bm{\zeta}_k|\mathcal{F}_k]=0,
~ \mathbb{E}[\|  \zeta_{i,k}\|^2|\mathcal{F}_k]\leq \nu_i^2,{\rm~and~}
\\&
\mathbb{E}[\| \bm{\zeta}_k\|^2|\mathcal{F}_k]\leq \sum_{i=1}^N \nu_i^2\triangleq \nu^2.
\end{split}
\end{equation}

\subsection{Main Results}

Define
\begin{align}
 \widehat{D}&(\mathbf{x},\mathbf{y})=\sum_{i=1}^N D(x_i,y_i)  ~{\rm for~any~} \mathbf{x},\mathbf{y}\in X,  \label{def-hatD3}
 \end{align}
With the definition of the Bregman's distance, we can replace the  strong monotonicity assumption  by the following assumption. This assumption  taken from \cite{kotsalis2020simple}  includes
 $\langle \phi(\mathbf{x}),  \mathbf{x}-\mathbf{x}^* \rangle \geq \mu \| \mathbf{x} -\mathbf{x}^*\|, \forall  \mathbf{x}  \in X $ as the special case when $h(x)=\|x\|^2/2$. 
\begin{assumption} \label{ass-monotone}
	There exists a constat $\mu>0$ such that
	\[\langle \phi(\mathbf{x}),  \mathbf{x}-\mathbf{x}^* \rangle \geq 2\mu \widehat{D}(\mathbf{x} ,\mathbf{x}^*) ,\quad \forall  \mathbf{x}  \in X   .\]
\end{assumption}
With this condition, we now  state a  convergence property of Algorithm \ref{alg1}, for which the proof can be found in Section \ref{subsec-prp1}.
\begin{proposition}\label{prp1} Consider Algorithm \ref{alg1}.  Let Assumptions \ref{ass-graph}-\ref{ass-monotone}  hold. Assume, in addition, that there exists a
positive sequence $\{\theta_k\}_{k\geq 1}$ satisfying
	\begin{align}
	& \theta_{k+1} \alpha_{k+1} \lambda_{k+1}=\theta_k \alpha_k ,\label{rela-theta}
	\\& \theta_{k-1}\geq 16   \alpha_k^2 \lambda_k^2  L^2\theta_k,\label{rela-theta2}
\\& \theta_k\leq \theta_{k-1}(2\mu \alpha_{k-1}+1) {\rm~and~}  8L^2\alpha_k^2\leq 1.\label{ass-theta}
\end{align}
	Then
	\begin{equation}\label{lemf5}
	\begin{split}
	&\theta_t (2\mu \alpha_t+1/2)\mathbb{E} [\widehat{D}(\mathbf{x}_{t+1},\mathbf{x}^*)]  +\sum_{k=1}^{t-1} {\theta_k \over 4} \mathbb{E} [\widehat{D}(\mathbf{x}_k,\mathbf{x}_{k+1}) ]\\ & \leq   \theta_1 \mathbb{E} [\widehat{D}(\mathbf{x}_1,\mathbf{x}^*)]+8\nu^2 \sum_{k=1}^t\theta_k \alpha_k^2 \lambda_k^2  +2   \theta_t\alpha_t^2 \nu^2
	\\&	  +    \theta_1 \alpha_1 \lambda_1    \mathbb{E} [  \varepsilon_1]
	+ 2 \sum_{k=1}^t \theta_k \alpha_k     \mathbb{E} [  \varepsilon_{k+1}],
	\end{split} \end{equation}
where
	\begin{equation}\label{def-vek}
	\begin{split}
	& \varepsilon_k\triangleq 2\sum_{i=1}^N L_{f_i} M_i \big( \| N \hat{v}_{i,k+1}-\sigma(\mathbf{x}_k)\|
	+\| N \hat{v}_{i,k }-\sigma(\mathbf{x}_{k-1})\| \big)  .
	\end{split} \end{equation}
\end{proposition}

\begin{remark} Consider the special case where the  digraph $\mathcal{G}_k$ is a complete graph for each time $k$ with $W_k={\mathbf{1}_N \mathbf{1}_N ^T \over N}$. Then \eqref{lemf5} becomes
	\begin{equation*}
	\begin{split}
	&\theta_t (2\mu \alpha_t+1/2)\mathbb{E} [\widehat{D}(\mathbf{x}_{t+1},\mathbf{x}^*)]  +\sum_{k=1}^{t-1} {\theta_k \over 4} \mathbb{E} [\widehat{D}(\mathbf{x}_k,\mathbf{x}_{k+1}) ]
	\\&\leq   \theta_1 \mathbb{E} [\widehat{D}(\mathbf{x}_1,\mathbf{x}^*)]+8\nu^2 \sum_{k=1}^t\theta_k \alpha_k^2 \lambda_k^2  +2   \theta_t\alpha_t^2 \nu^2 .
	\end{split} \end{equation*}
	This recovers the bound of   \cite[Theorem 3.3]{kotsalis2020simple}.
\end{remark}

In the following, we  establish a  bound on the consensus error $\|\sigma(\mathbf{x}_k) -N\hat{v}_{i,k+1}\|$ of the aggregate,
for which the proof can be found in Section \ref{subsec-prp2}.
\begin{proposition} \label{prp2} Consider Algorithm \ref{alg1}. Let Assumptions \ref{ass-graph}, \ref{ass-payoff}, and \ref{ass-hF}  hold.  Then
	\begin{equation}\label{agg-bd0}
	\begin{split}
	&\|\sigma(\mathbf{x}_k) -N\hat{v}_{i,k+1}\|\leq \theta M_HN \beta^{k}  \\
	&  +\theta N \sum_{s=1}^k \beta^{k-s}   \alpha_{s-1}  \sum_{i=1}^N  \Big((1+2\lambda_{s-1})C_i +\lambda_{s-1}\| \zeta_{i,s-2}\|\\
	&\qquad\qquad\qquad\qquad \quad + (1+\lambda_{s-1})\| \zeta_{i,s-1}\|\Big),
	\end{split}
	\end{equation}
	where the constants $\theta, \beta$ are defined in \eqref{geometric}, and
\begin{align}C_i  \triangleq  N\tilde{C}   L_{fi}+\max_{  \mathbf{x}\in X}\|\phi_i( \mathbf{x})\| .\label{def-C}
\end{align}
\end{proposition}

By combining Proposition  \ref{prp1} and Proposition  \ref{prp2},
we can  show that the proposed  method   can achieve the optimal convergence rate for solving the stochastic smooth and strongly monotone aggregative games.
The proof can be found in Section \ref{subsec-thm1}.
\begin{theorem}\label{thm1} Consider Algorithm \ref{alg1}. Suppose Assumptions \ref{ass-graph}-\ref{ass-monotone} hold. Set
	\begin{align*}
	&c_0={4L\over \mu},~\alpha_k={1\over \mu(k+c_0-1)},
	\\& \theta_k=(k+c_0+1)(k+c_0 ), ~\lambda_k
	={\theta_{k-1}\alpha_{k-1} \over \theta_k\alpha_k}.
	\end{align*}
	Then the following hold with $c_e  \triangleq  4N\tilde{C} \sum_{i=1}^N L_{f_i} M_i$.
	\begin{equation}\label{result3}
	\begin{split}
	&\mathbb{E} [\widehat{D}(\mathbf{x}_{t+1},\mathbf{x}^*)]
	\leq   {2 (c_0+2)(c_0+1)  \widehat{D}(\mathbf{x}_1,\mathbf{x}^*)  \over (t+c_0+1)(t+c_0 )}
	 \\&+{8( \nu^2 + c_1)\over \mu^2(t+c_0+1)(t+c_0 )}
	+ {2c_0 (c_0+1)c_e  \over \mu (c_0-1) (t+c_0+1)(t+c_0 )}
	\\&+ {8( c_2+4\nu^2)t \over \mu^2(t+c_0+1)(t+c_0 ) }    ,
	\end{split} \end{equation}
where \begin{align} c_1 &\triangleq  2\theta M_HN (1+\beta)\mu\beta \left({  (c_0-1)  \over   1- \beta  }+
{1\over ( 1- \beta )^2 }\right)  \sum_{i=1}^N L_{f_i} M_i,\label{def-C1}
\\      c_2 &\triangleq   {8\theta N (3-\beta)\sum_{i=1}^N   (C_i + \nu_i)    \over  ( 1- \beta )^2 } \sum_{i=1}^N L_{f_i} M_i.\label{def-C2}
\end{align}
\end{theorem}

\begin{cor} \label{cor1}
 The number of iterations (the same as communication rounds) required by Algorithm \ref{alg1} for obtaining an approximate Nash equilibrium $\bar{\mathbf{x}} \in X$ satisfying $  \mathbb{E} [\widehat{D}(\bar{\mathbf{x}},\mathbf{x}^*)]\leq \epsilon$ is bounded by
	\[ \max \left( {L \sqrt{\widehat{D}(\mathbf{x}_1,\mathbf{x}^*)} \over \mu \sqrt{\epsilon}},
	{\sqrt{ \nu^2 + c_1 }\over \mu \sqrt{\epsilon}},{\sqrt{L} \over \mu \sqrt{\epsilon}}, {  c_2+4\nu^2  \over \mu^2 \epsilon}
	\right).\]
\end{cor}

\section{Proof of Main Results}\label{sec:proof}

\subsection{Preliminary Results}


We now state a property from \cite[Proposition B.3]{bravo2018bandit},
and a well-known technical results regarding  the optimality condition of \eqref{alg-strategy0} (see   \cite[Lemma 3.1]{lan2020first}).
\begin{lemma}\label{lemr}
	Let $h$ be a smooth and 1-strongly convex regularizer.
	Then
	\begin{align}\label{property-D}
	D (x_i,y_i)\geq \frac{1 }{2}\|x_i - y_i\|^2  , \quad \forall x_i,y_i\in X_i,  \end{align}
	and \begin{align}
	& \alpha_k \big \langle
	(1+ \lambda_k ) q_i(x_{i,k}, N\hat{v}_{i,k+1};\xi_{i,k}),x_{i,k+1}-x_i  \big  \rangle \notag
	\\&- \alpha_k \big \langle \lambda_k q_i(x_{i,k-1}, N\hat{v}_{i,k };\xi_{i,k-1})  ,x_{i,k+1}-x_i  \big  \rangle + D(x_{i,k},x_{i,k+1}) \notag
	\\&\leq   D (x_{i,k},x_i)- D(x_{i,k+1},x_i),
	\quad \forall x_i\in X_i. \label{low-bdd}
	\end{align}
\end{lemma}

Then from \eqref{def-hatD3} and \eqref{property-D} it follows that
\begin{align}\label{property-Dhat}
\widehat{D}(\mathbf{x},\mathbf{y}) \geq {1\over 2} \|\mathbf{x}-\mathbf{y}\|^2,
\quad \forall  \mathbf{x},\mathbf{y}\in X.
\end{align}
Furthermore, we define
\begin{align}
& \Delta D_{i,k}(x_i)\triangleq D (x_{i,k},x_i)- D(x_{i,k+1},x_i),  \label{def-hatD1}
\\& \Delta q_{i,k} \triangleq  q_i(x_{i,k}, N\hat{v}_{i,k+1};\xi_{i,k}) -   q_i(x_{i,k-1}, N\hat{v}_{i,k };\xi_{i,k-1})   ,  \label{def-hatD2}
\\& \Delta \widehat{D}_k(\mathbf{x})= \sum_{i=1}^N \Delta D_{i,k}(x_i)=
\widehat{D}(\mathbf{x}_k,\mathbf{x})- \widehat{D}(\mathbf{x}_{k+1},\mathbf{x}). \label{def-hatD4}
\end{align}

We are now ready to show the convergence properties  of the proposed  method.

\begin{lemma}
	Let    $\{\mathbf{x}_k\}$ be generated by Algorithm \ref{alg1}.
Let Assumptions   \ref{ass-payoff} and \ref{ass-hF} hold. Suppose, in addition,  that \eqref{rela-theta}  and \eqref{rela-theta2} hold for some positive sequence $\{\theta_k\}_{k\geq 1}$.

	Then
	\begin{equation}\label{lemf}
	\begin{split}
	&\sum_{k=1}^t \theta_k \left( \alpha_k  \sum_{i=1}^N\big \langle
	q_i(x_{i,k+1}, N\hat{v}_{i,k+2};\xi_{i,k+1}), x_{i,k+1}-x_i \big  \rangle
	+\right.\\&\left.\widehat{D}(\mathbf{x}_{k+1},\mathbf{x}) +  {1\over 4} \widehat{D}(\mathbf{x}_k,\mathbf{x}_{k+1})  \right)
	\\ &    - 2\alpha_t^2\theta_t L^2  \| \mathbf{x}_{ t+1 }-\mathbf{x}  \|^2
	-   \alpha_t \theta_t \langle \bm{\zeta}_{ t+1}-     \bm{\zeta}_t , \mathbf{x}_{ t+1 }-\mathbf{x} \big  \rangle  \\& \leq \sum_{k=1}^t \theta_k \widehat{D}(\mathbf{x}_k,\mathbf{x})
	+  \sum_{k=1}^t \theta_k \alpha_k \lambda_k      \varepsilon_k
	\\&+2\sum_{k=1}^t\theta_k \alpha_k^2 \lambda_k^2 \| \bm{\zeta}_k-\bm{\zeta}_{ k-1}\|^2
	+ \alpha_t \theta_t \varepsilon_{t+1} .
	\end{split}
	\end{equation}
\end{lemma}

{\bf Proof.}  With the definitions \eqref{def-hatD1} and \eqref{def-hatD2},
\eqref{low-bdd} becomes
\begin{equation}\label{low-bdd2}
\begin{split}
\Delta D_{i,k}(x_i)& \geq \alpha_k \big \langle
q_i(x_{i,k+1}, N\hat{v}_{i,k+2};\xi_{i,k+1}), x_{i,k+1}-x_i \big  \rangle
\\&-\alpha_k \big \langle \Delta q_{i,k+1}  , x_{i,k+1}-x_i \big  \rangle
\\& +\alpha_k \lambda_k   \langle \Delta q_{i,k}  , x_{i,k }-x_i \big  \rangle
\\&+\alpha_k \lambda_k   \langle \Delta q_{i,k}  , x_{i,k+1 }-x_{i,k } \big  \rangle + D(x_{i,k},x_{i,k+1})  .
\end{split} \end{equation}

By multiplying both sides of the above inequality by $\theta_k$, and summing up from
$i=1$ to $N$, and   $k=1$ to $t $,  we obtain that
\begin{align} \label{low-bdd3}
&\sum_{k=1}^t \theta_k \Delta \widehat{D}_k(\mathbf{x}) \geq \sum_{k=1}^t \theta_k\alpha_k \lambda_k  \sum_{i=1}^N \langle \Delta q_{i,k}  , x_{i,k }-x_i \big  \rangle \notag
\\& + \sum_{k=1}^t \theta_k \alpha_k \sum_{i=1}^N \big \langle \notag
q_i(x_{i,k+1}, N\hat{v}_{i,k+2};\xi_{i,k+1}), x_{i,k+1}-x_i \big  \rangle  \notag
\\& -\sum_{k=1}^t \theta_k \alpha_k  \sum_{i=1}^N \big \langle \Delta q_{i,k+1}  , x_{i,k+1}-x_i \big  \rangle
\\& + \sum_{k=1}^t \theta_k \sum_{i=1}^N D(x_{i,k},x_{i,k+1})+\sum_{k=1}^t \theta_k \alpha_k \lambda_k  \sum_{i=1}^N \langle \Delta q_{i,k}  , x_{i,k+1 }-x_{i,k } \big  \rangle .\notag
\end{align}

Then by recalling that  $\Delta q_{i,1}=0$,  using   \eqref{def-hatD3} and \eqref{rela-theta},
we obtain
\begin{align}
& \sum_{k=1}^t \theta_k \Delta \widehat{D}_k(\mathbf{x}) \geq  -  \theta_t \alpha_t  \sum_{i=1}^N\big  \langle \Delta q_{i,t+1}  , x_{i,t+1}-x_i \big  \rangle  \notag
\\& + \underbrace{ \sum_{k=1}^t \theta_k \widehat{D}(\mathbf{x}_k,\mathbf{x}_{k+1}) +\sum_{k=1}^t \theta_k \alpha_k \lambda_k  \sum_{i=1}^N \langle \Delta q_{i,k}  , x_{i,k+1 }-x_{i,k } \big  \rangle  }_{\triangleq Q_t}  \label{low-bdd4}
\\&+ \sum_{k=1}^t \theta_k \alpha_k
\sum_{i=1}^N\big \langle
q_i(x_{i,k+1}, N\hat{v}_{i,k+2};\xi_{i,k+1}), x_{i,k+1}-x_i \big  \rangle  \notag. \notag
\end{align}

By the definitions of \eqref{def-zeta} and \eqref{def-hatD2}, we have
\begin{equation}\label{bd1}
\begin{split}
& \sum_{i=1}^N  \langle \Delta q_{i,k}  , x_{i,k+1 }-x_{i,k } \big  \rangle=
\sum_{i=1}^N  \langle  \zeta_{i,k}-     \zeta_{i,k-1}  , x_{i,k+1 }-x_{i,k } \big  \rangle \\&+\sum_{i=1}^N
\langle  F_i(x_{i,k}, N\hat{v}_{i,k+1}) -   F_i(x_{i,k-1}, N\hat{v}_{i,k } )    , x_{i,k+1 }-x_{i,k } \big  \rangle
\\& =   \langle  \bm{\zeta}_{ k}-     \bm{\zeta}_{ k-1}  , \mathbf{x}_{ k+1 }-\mathbf{x}_k \big  \rangle
\\&+\sum_{i=1}^N \langle  F_i(x_{i,k}, N\hat{v}_{i,k+1})-  F_i(x_{i,k-1}, N\hat{v}_{i,k } )    , x_{i,k+1 }-x_{i,k } \big  \rangle.
\end{split} \end{equation}

By using   \eqref{def-bdst} and Assumption \ref{ass-hF}, we derive
\begin{equation}\label{bd-fd}
\begin{split}
& \sum_{i=1}^N \langle  F_i(x_{i,k}, N\hat{v}_{i,k+1}) -   F_i(x_{i,k-1}, N\hat{v}_{i,k } )   , x_{i,k+1 }-x_{i,k }  \rangle
\\   & \overset{\eqref{def-phi}}{=}  \sum_{i=1}^N\langle F_i(x_{i,k}, N\hat{v}_{i,k+1}) -F_i(x_{i,k}, \sigma(\mathbf{x}_k))   , x_{i,k+1 }-x_{i,k }  \rangle
\\&+ \sum_{i=1}^N \langle  F_i(x_{i,k-1}, \sigma(\mathbf{x}_{k-1})) - F_i(x_{i,k-1}, N\hat{v}_{i,k } )  , x_{i,k+1 }-x_{i,k }  \rangle
\\&+ \sum_{i=1}^N \langle \phi_i(\mathbf{x}_k)-\phi_i(\mathbf{x}_{k-1}), x_{i,k+1 }-x_{i,k }  \rangle
\\ &\geq - \sum_{i=1}^N \| F_i(x_{i,k}, N\hat{v}_{i,k+1}) -F_i(x_{i,k}, \sigma(\mathbf{x}_k))\| \| x_{i,k+1 }-x_{i,k } \| \\&
- \sum_{i=1}^N  \| F_i(x_{i,k-1}, \sigma(\mathbf{x}_{k-1})) - F_i(x_{i,k-1}, N\hat{v}_{i,k } ) \| \| x_{i,k+1 }-x_{i,k }  \|
\\&+   \langle \phi (\mathbf{x}_k)-\phi (\mathbf{x}_{k-1}), \mathbf{x}_{ k+1 }-\mathbf{x}_{ k }  \rangle
\\& \geq  -  2\sum_{i=1}^N L_{f_i} M_i \big( \| N \hat{v}_{i,k+1}-\sigma(\mathbf{x}_k)\|
	+\| N \hat{v}_{i,k }-\sigma(\mathbf{x}_{k-1})\| \big)
\\& -L \|\mathbf{x}_k-\mathbf{x}_{k-1} \| \| \mathbf{x}_{ k+1 }-\mathbf{x}_{ k }  \|.
\end{split} \end{equation}
This together with \eqref{bd1}, and the definitions of $Q_t$ and  $ \varepsilon_k $ in \eqref{low-bdd4} and \eqref{def-vek}
implies that
\begin{equation}
\begin{split}
Q_t&\geq  \sum_{k=1}^t \theta_k \widehat{D}(\mathbf{x}_k,\mathbf{x}_{k+1}) +\sum_{k=1}^t \theta_k \alpha_k \lambda_k \langle  \bm{\zeta}_{ k}-     \bm{\zeta}_{ k-1}  , \mathbf{x}_{ k+1 }-\mathbf{x}_{ k }  \big  \rangle
\\&- \sum_{k=1}^t \theta_k \alpha_k \lambda_k  L \|\mathbf{x}_k-\mathbf{x}_{k-1} \| \| \mathbf{x}_{ k+1 }-\mathbf{x}_{ k }  \|
- \sum_{k=1}^t \theta_k \alpha_k \lambda_k      \varepsilon_k .
\end{split} \end{equation}
Then by  \eqref{property-Dhat}, the following holds  with defining $\theta_0\triangleq 0.$
\begin{align}\label{bd2}
Q_t&\geq  {1\over 4}\sum_{k=1}^t \theta_k \widehat{D}(\mathbf{x}_k,\mathbf{x}_{k+1})  \notag
\\&+\sum_{k=1}^t
\left( \underbrace{\theta_k \alpha_k \lambda_k \langle \bm{\zeta}_{ k}-     \bm{\zeta}_{ k-1}  , \mathbf{x}_{ k+1 }-\mathbf{x}_{ k } \big  \rangle + {\theta_k \over 8} \| \mathbf{x}_{ k+1 }-\mathbf{x}_{ k }  \|^2}_{\triangleq {\rm ~Term}~1} \right)\notag
\\&+\sum_{k=1}^t \left( - \theta_k \alpha_k \lambda_k  L \|\mathbf{x}_k-\mathbf{x}_{k-1} \| \| \mathbf{x}_{ k+1 }-\mathbf{x}_{ k }  \|\right) \notag \\&+\sum_{k=1}^t\left( {\theta_k \over 8} \| \mathbf{x}_{ k+1 }-\mathbf{x}_{ k }  \|^2 +{\theta_{k-1} \over 8} \| \mathbf{x}_{k-1}-\mathbf{x}_{ k }  \|^2  \right) \notag
\\&+{\theta_t \over 8} \| \mathbf{x}_{ t+1 }-\mathbf{x}_t  \|^2
-  \sum_{k=1}^t \theta_k \alpha_k \lambda_k \varepsilon_k .
\end{align}

We let
\begin{align}
{\rm Term}~2 \triangleq & - \theta_k \alpha_k \lambda_k  L \|\mathbf{x}_k-\mathbf{x}_{k-1} \| \| \mathbf{x}_{ k+1 }-\mathbf{x}_{ k }  \|\ \notag
\\&+ {\theta_k \over 8} \| \mathbf{x}_{ k+1 }-\mathbf{x}_{ k }  \|^2 +{\theta_{k-1} \over 8} \| \mathbf{x}_{k-1}-\mathbf{x}_{ k }  \|^2.   \notag
\end{align}
By using \eqref{rela-theta2} and   the Young's inequality $a^2+b^2\geq 2\|a\|\|b\|$, we derive ${\rm Term}~2  \geq 0 $. Also, we obtain
${\rm ~Term}~1  \geq -2\theta_k \alpha_k^2 \lambda_k^2 \|\bm{\zeta}_{ k}-     \bm{\zeta}_{ k-1}  \|^2 $. Then by substituting the bounds of Term 1 and Term 2 into \eqref{bd2}, we have
\begin{equation}\label{bd3}
\begin{split}
Q_t\geq & {1\over 4}\sum_{k=1}^t \theta_k \widehat{D}(\mathbf{x}_k,\mathbf{x}_{k+1})
-2\sum_{k=1}^t\theta_k \alpha_k^2 \lambda_k^2 \|\bm{\zeta}_{ k}-     \bm{\zeta}_{ k-1}  \|^2
\\&+{\theta_t \over 8} \| \mathbf{x}_{ t+1 }-\mathbf{x}_t  \|^2
-  \sum_{k=1}^t \theta_k \alpha_k \lambda_k      \varepsilon_k .
\end{split} \end{equation}
This together with  \eqref{low-bdd4} implies that
\begin{equation}\label{bd4}
\begin{split}
& \sum_{k=1}^t \theta_k \Delta \widehat{D}_k(\mathbf{x}) \geq {1\over 4}\sum_{k=1}^t \theta_k \widehat{D}(\mathbf{x}_k,\mathbf{x}_{k+1})
-2\sum_{k=1}^t\theta_k \alpha_k^2 \lambda_k^2 \|\bm{\zeta}_{ k}-     \bm{\zeta}_{ k-1}  \|^2
\\& +  \sum_{k=1}^t \theta_k \alpha_k
\sum_{i=1}^N\big \langle
q_i(x_{i,k+1}, N\hat{v}_{i,k+2};\xi_{i,k+1}), x_{i,k+1}-x_i \big  \rangle
\\& -  \sum_{k=1}^t \theta_k \alpha_k \lambda_k      \varepsilon_k
+{\theta_t \over 8} \| \mathbf{x}_{ t+1 }-\mathbf{x}_t  \|^2 -  \theta_t \alpha_t  \sum_{i=1}^N\big  \langle \Delta q_{i,t+1}  , x_{i,t+1}-x_i \big  \rangle .
\end{split} \end{equation}

Similarly to \eqref{bd1}, we have that
\begin{equation}
\begin{split}
& -\sum_{i=1}^N  \langle \Delta q_{i,t+1}  , x_{i,t+1 }-x_{i  } \big  \rangle=
-  \langle  \bm{\zeta}_{ t+1}-     \bm{\zeta}_t  , \mathbf{x}_{ t+1 }-\mathbf{x} \big  \rangle
\\& -\sum_{i=1}^N    \langle  F_i(x_{i,t+1}, N\hat{v}_{i,t+2}) -   F_i(x_{i,t }, N\hat{v}_{i,t+1 } )    , x_{i,t+1 }-x_i \big  \rangle
\end{split} \end{equation}
Similarly to \eqref{bd-fd}, by using  Assumption \ref{ass-hF},
\eqref{def-phi}, and \eqref{def-bdst}, we have that
\begin{equation*}
\begin{split}
& -\sum_{i=1}^N \langle  F_i(x_{i,t+1}, N\hat{v}_{i,t+2}) -   F_i(x_{i,t }, N\hat{v}_{i,t+1 } )   , x_{i,t+1 }-x_i  \rangle
\\  & \geq  -    \varepsilon_{t+1}
-L \|\mathbf{x}_{t+1}-\mathbf{x}_t \| \| \mathbf{x}_{ t+1 }-\mathbf{x}  \|.
\end{split} \end{equation*}
Therefore,
\begin{equation*}
\begin{split}
&{1 \over 8} \| \mathbf{x}_{ t+1 }-\mathbf{x}_t  \|^2 -   \alpha_t  \sum_{i=1}^N\big  \langle \Delta q_{i,t+1}  , x_{i,t+1}-x_i \big  \rangle
\\& \geq {1 \over 8} \| \mathbf{x}_{ t+1 }-\mathbf{x}_t  \|^2 -  \alpha_t \langle \bm{\zeta}_{ t+1}-     \bm{\zeta}_t , \mathbf{x}_{ t+1 }-\mathbf{x} \big  \rangle
\\&-   \alpha_t \varepsilon_{t+1}
- \alpha_t L \|\mathbf{x}_{t+1}-\mathbf{x}_t \| \| \mathbf{x}_{ t+1 }-\mathbf{x}  \|
\\& \geq -   \alpha_t \langle  \bm{\zeta}_{ t+1}-     \bm{\zeta}_t , \mathbf{x}_{ t+1 }-\mathbf{x} \big  \rangle  - 2\alpha_t^2 L^2  \| \mathbf{x}_{ t+1 }-\mathbf{x}  \|^2-    \alpha_t \varepsilon_{t+1},
\end{split}
\end{equation*}
where the last inequality follows by the Young' s inequality.
This incorporating with \eqref{bd4} produces
\begin{equation*}
\begin{split}
& \sum_{k=1}^t \theta_k \Delta \widehat{D}_k(\mathbf{x}) \geq {1\over 4}\sum_{k=1}^t \theta_k \widehat{D}(\mathbf{x}_k,\mathbf{x}_{k+1})
-2\sum_{k=1}^t\theta_k \alpha_k^2 \lambda_k^2 \| \bm{\zeta}_{ k}-     \bm{\zeta}_{ k-1}\|^2
\\& + \sum_{k=1}^t \theta_k \alpha_k
\sum_{i=1}^N\big \langle
q_i(x_{i,k+1}, N\hat{v}_{i,k+2};\xi_{i,k+1}), x_{i,k+1}-x_i \big  \rangle
\\& - \sum_{k=1}^t \theta_k \alpha_k \lambda_k      \varepsilon_k
-   \alpha_t \theta_t\langle  \bm{\zeta}_{ k+1}-     \bm{\zeta}_{ k  }  , \mathbf{x}_{ t+1 }-\mathbf{x} \big  \rangle
\\& - 2\alpha_t^2\theta_t L^2  \| \mathbf{x}_{ t+1 }-\mathbf{x}  \|^2
-    \alpha_t\theta_t \varepsilon_{t+1} .
\end{split} \end{equation*}
By recalling the definition  \eqref{def-hatD4}  and rearranging the terms,
we  prove  the lemma.
\hfill $\Box$

\subsection{Proof of Proposition \ref{prp1}}\label{subsec-prp1}

  Since $\mathbf{x}_{k }$ is adapted to  $\mathcal{F}_k,$ from  \eqref{bd-noise}it follows  that
  \[\mathbb{E}[\langle \bm{\zeta}_{ k  }  , \mathbf{x}_{ k  }-\mathbf{x} \big  \rangle ]
=\mathbb{E}\big[ \mathbb{E} [\langle \bm{\zeta}_{ k  }  , \mathbf{x}_{ k  }-\mathbf{x} \big  \rangle | \mathcal{F}_k]\big] =\mathbb{E}\big[ \mathbb{E} [\langle \bm{\zeta}_{ k  } | \mathcal{F}_k] , \mathbf{x}_{ k  }-\mathbf{x} \big  \rangle \big]=0 \] for any $k\geq 1.$ Therefore,
\begin{align}\label{re-zeta}
 &\mathbb{E}[\langle  \bm{\zeta}_{ t+1}-     \bm{\zeta}_t  , \mathbf{x}_{ t+1 }-\mathbf{x} \big  \rangle ]=\mathbb{E}[\langle   - \bm{\zeta}_t , \mathbf{x}_{ t+1 }-\mathbf{x} \big  \rangle ] \notag
\\&=\mathbb{E}[\langle-\bm{\zeta}_t  , \mathbf{x}_{ t+1}-\mathbf{x}_t \big  \rangle ].
\end{align}
By noting that $ q_i(x_{i,k+1}, N\hat{v}_{i,k+2};\xi_{i,k+1}) =F_i(x_{i,k+1}, N\hat{v}_{i,k+2})+\zeta_{i,k+1},$ we have
\begin{align}\label{re-qik} &
\mathbb{E}[\langle q_i(x_{i,k+1}, N\hat{v}_{i,k+2};\xi_{i,k+1}), x_{i,k+1}-x_i \big  \rangle ] \notag
\\&= \mathbb{E}[\langle
F_i(x_{i,k+1}, N\hat{v}_{i,k+2} ), x_{i,k+1}-x_i \big  \rangle ].
\end{align}
In addition, by \eqref{bd-noise} we have
\[\mathbb{E} \big[\| \bm{\zeta}_{ k}-     \bm{\zeta}_{ k-1}\|^2 \big]
\leq 2\mathbb{E} \big[\| \bm{\zeta}_{ k} \|^2 \big] +2\mathbb{E} \big
[\|\bm{\zeta}_{ k-1}\|^2 \big]\leq 4\nu^2 .\]

Then by taking unconditional expectations on both sides of \eqref{lemf}, using \eqref{re-zeta} and \eqref{re-qik}, we obtain that
 \begin{equation}\label{lemf2}
\begin{split}
&\sum_{k=1}^t \theta_k \left( \alpha_k
\sum_{i=1}^N\mathbb{E}[\langle
F_i(x_{i,k+1}, N\hat{v}_{i,k+2} ), x_{i,k+1}-x_i \big  \rangle ]\right)
\\&+\sum_{k=1}^t \theta_k \left(\mathbb{E} [\widehat{D}(\mathbf{x}_{k+1},\mathbf{x})] +  {1\over 4} \mathbb{E} [\widehat{D}(\mathbf{x}_k,\mathbf{x}_{k+1}) ] \right)
\\ &    - 2\alpha_t^2\theta_t L^2  \mathbb{E} [\| \mathbf{x}_{ t+1 }-\mathbf{x}  \|^2]
+ \alpha_t \theta_t \mathbb{E}[\langle  \bm{\zeta}_{ t  }  , \mathbf{x}_{ t+1 }-\mathbf{x}_t \big  \rangle ]  \\& \leq \sum_{k=1}^t \theta_k \mathbb{E} [\widehat{D}(\mathbf{x}_k,x)]
+   \sum_{k=1}^t \theta_k \alpha_k \lambda_k    \mathbb{E} [  \varepsilon_k]
\\&+8\nu^2\sum_{k=1}^t\theta_k \alpha_k^2 \lambda_k^2  +\alpha_t\theta_t \mathbb{E} [\varepsilon_{t+1}] .
\end{split} \end{equation}

Note by   Assumption \ref{ass-hF}(b), \eqref{def-phi}, and \eqref{def-bdst} that
\begin{equation*}
\begin{split}
& \langle   F_i(x_{i,k}, N\hat{v}_{i,k+1})  , x_{i,k}- x_i \rangle
\\&  =  \langle   F_i(x_{i,k},
\sigma(\mathbf{x}_k)), x_{i,k}- x_i \rangle
\\&+\langle    F_i(x_{i,k}, N\hat{v}_{i,k+1}) -  F_i(x_{i,k},\sigma(\mathbf{x}_k)), x_{i,k}- x_i \rangle
\\ & \geq  \langle\phi_i( \mathbf{x}_k), x_{i,k}- x_i \rangle
\\&-\| F_i(x_{i,k}, \sigma(\mathbf{x}_k)) -F_i(x_{i,k}, N\hat{v}_{i,k+1})\|
\| x_{i,k}- x_i\|
\\ & \geq \langle\phi_i( \mathbf{x}_k), x_{i,k}- x_i \rangle-2L_{fi}M_i\| \sigma(\mathbf{x}_k)- N\hat{v}_{i,k+1}) \| .
\end{split} \end{equation*}
Then by recalling the definition of $\varepsilon_k$ in  \eqref{def-vek},
and using Assumption \ref{ass-monotone}, we obtain that
\begin{equation}\label{bd-fx}
\begin{split}
 \sum_{i=1}^N \langle  F_i(x_{i,k}, N\hat{v}_{i,k+1})  , x_{i,k}- x_i^* \rangle&
\geq  \langle\phi( \mathbf{x}_k), \mathbf{x}_k- \mathbf{x}^*  \rangle-2\varepsilon_k
\\&\geq 2\mu \widehat{D}(\mathbf{x}_k,\mathbf{x}^*) - \varepsilon_k.
\end{split} \end{equation}

Note by \eqref{property-Dhat} and \eqref{bd-noise} that
\begin{equation*}
\begin{split}
&  {1\over 4} \mathbb{E} [\widehat{D}(\mathbf{x}_t,\mathbf{x}_{t+1}) ]  +\alpha_t  \mathbb{E}[\langle  \bm{\zeta}_{ t  }  , \mathbf{x}_{ t+1 }-\mathbf{x}_k \big  \rangle ]
\\& \geq {1\over 8}\mathbb{E} [\| \mathbf{x}_t-\mathbf{x}_{t+1} \|^2] +\alpha_t  \mathbb{E}[\langle  \bm{\zeta}_{ t }  , \mathbf{x}_{ t+1 }-\mathbf{x}_t \big  \rangle ]
\\& \geq {1\over 8}\mathbb{E} [\| \mathbf{x}_t-\mathbf{x}_{t+1} \|^2] -\alpha_t
\sqrt{\mathbb{E} [\| \mathbf{x}_t-\mathbf{x}_{t+1} \|^2]}
\sqrt{\mathbb{E}[\| \bm{\zeta}_{ t  } \|^2    ]}
\\& \geq
-2\alpha_t^2 \mathbb{E}[\| \bm{\zeta}_{ t  } \|^2    ] \geq -2\alpha_t^2 \nu^2.
\end{split} \end{equation*}
This together with \eqref{lemf2} and \eqref{bd-fx} implies that
\begin{equation}\label{lemf3}
\begin{split}
&\sum_{k=1}^t \theta_k   (2\mu \alpha_k+1)\mathbb{E} [\widehat{D}(\mathbf{x}_{k+1},\mathbf{x}^*)]
\\&+\sum_{k=1}^{t-1} {\theta_k \over 4} \mathbb{E} [\widehat{D}(\mathbf{x}_k,\mathbf{x}_{k+1}) ] - 2\alpha_t^2\theta_t L^2  \mathbb{E} [\| \mathbf{x}_{ t+1 }-\mathbf{x}^*  \|^2]
\\ & \leq \sum_{k=1}^t \theta_k \mathbb{E} [\widehat{D}(\mathbf{x}_k,\mathbf{x}^*)]
+   \sum_{k=1}^t \theta_k \alpha_k \lambda_k    \mathbb{E} [  \varepsilon_k]+2  \theta_t\alpha_t^2 \nu^2
\\&+8\nu^2 \sum_{k=1}^t\theta_k \alpha_k^2 \lambda_k^2     +  \alpha_t \theta_t \mathbb{E} [\varepsilon_{t+1}]+  {\sum_{k=1}^t \theta_k \alpha_k     \mathbb{E} [  \varepsilon_{k+1}]}.
\end{split} \end{equation}

By recalling   $\theta_k\leq \theta_{k-1}(2\mu \alpha_{k-1}+1) $ from the  condition \eqref{ass-theta}, we obtain   that
\begin{equation}\label{lemf4}
\begin{split}
&\sum_{k=1}^t \theta_k   (2\mu \alpha_k+1)\mathbb{E} [\widehat{D}(\mathbf{x}_{k+1},\mathbf{x}^*)]   -  2\alpha_t^2\theta_t L^2  \mathbb{E} [\| \mathbf{x}_{ t+1 }-\mathbf{x}^*  \|^2]
\\& \geq \theta_t \left( (2\mu \alpha_t+1/2)\mathbb{E} [\widehat{D}(\mathbf{x}_{t+1},\mathbf{x}^*)]  \right)+
\sum_{k=1}^{t-1} \theta_{k+1} \mathbb{E} [\widehat{D}(\mathbf{x}_{k+1},\mathbf{x}^*)]
\\&+{\theta_t  \over 2}   \left( \mathbb{E} [\widehat{D}(\mathbf{x}_{t+1},\mathbf{x}^*)]  -   4\alpha_t^2  L^2  \mathbb{E} [\| \mathbf{x}_{ t+1 }-\mathbf{x}^*  \|^2]\right)
\\& \geq \theta_t  (2\mu \alpha_t+1/2)\mathbb{E} [\widehat{D}(\mathbf{x}_{t+1},\mathbf{x}^*)] + \sum_{k=1}^{t-1} \theta_{k+1} \mathbb{E} [\widehat{D}(\mathbf{x}_{k+1},\mathbf{x}^*)] ,
\end{split} \end{equation}
where the last inequality follows by \eqref{property-Dhat} and   $  8L^2\alpha_k^2\leq 1$ of \eqref{ass-theta}.  This together with \eqref{lemf3} and \eqref{rela-theta} proves the result.
\hfill $\Box$

\subsection{ Proof of Proposition \ref{prp2}.} \label{subsec-prp2}

Since $ v_{ i,1} = x_{i,1}  $ and $W(k)$ is doubly stochastic,
similarly to  \cite[Lemma 2]{koshal2016distributed},
we can show by induction that
\begin{align}\label{equi}  \sum_{i=1}^N v_{i,k} =\sum_{i=1}^N  x_{i,k} =\sigma(\mathbf{x}_k),\quad \forall k\geq 0.
\end{align}

Akin to   \cite[Eqn. (16)]{koshal2016distributed},
we obtain the following bound.
\begin{align*}
\left \|{ \sigma(\mathbf{x}_k)\over N} -\hat{v}_{i,k+1} \right \|\leq& \sum_{j=1}^n \left| {1\over N}-[\Phi(k,1)]_{ij}\right| \|v_{j,1}\| \\ &+\sum_{s=1}^k \sum_{j=1}^N
\left| {1\over N}-[\Phi(k,s)]_{ij}\right|  \big\|  x_{j,s}-x_{j,s-1}   \big \|.
\end{align*}
Then by using \eqref{geometric}  and  $v_{i,1}= x_{i,1} $,  we obtain that
\begin{equation}\label{bd-consensus1}
\begin{split}
\left \|{ \sigma(\mathbf{x}_k)\over N} -\hat{v}_{i,k+1} \right  \| \leq  \theta \beta^{k}\sum_{j=1}^N \| x_{j,1} \|   +   \theta \sum_{s=1}^k \beta^{k-s}   \sum_{j=1}^N   \big\|  x_{j,s}- x_{j,s-1}   \big \| .
\end{split} \end{equation}
This combined with \eqref{def-bdst} proves
\begin{align} \label{bd-con-err}
\left \|{ \sigma(\mathbf{x}_k)\over N} -\hat{v}_{i,k+1} \right\|
&\leq \theta \beta^{k}M_H  +  2 \theta \sum_{s=1}^k \beta^{k-s}   \sum_{j=1}^N   M_j \notag
\\&  \leq  \theta   M_H  +  2 \theta    \sum_{j=1}^N   M_j/(1-\beta)
\overset{\eqref{def-tC}} {=}\tilde{C}.
\end{align}

By  \eqref{def-bdst} and \eqref{equi}, we have that $\|\sigma(\mathbf{x}_k)\| \leq M_H$ for any $k\geq 0$.
Thus by \eqref{bd-con-err},  we obtain that  for each $i\in \mathcal{N} $,
\[\|N\hat{v}_{i,k+1}\| \leq N\tilde{C}+M_H,\quad \forall k\geq 0.\]
Then by using Assumption \ref{ass-hF}(b) and \eqref{def-phi},  we obtain that for   each $j\in \mathcal{N}$  and any $s\geq 0:$
\begin{align}\label{bf-F}
\| F_j(x_{j,s }, N\hat{v}_{j,s }) \|  &\leq \| F_j(x_{j,s }, N\hat{v}_{j,s }) -F_j(x_{j,s },\sigma(\mathbf{x}_s)\| \notag
\\&+\|F_j(x_{j,s },\sigma(\mathbf{x}_s)\| \notag
\\& \leq L_{fj} \|   N\hat{v}_{j,s }  - \sigma(\mathbf{x}_s)\| +\|\phi_j( \mathbf{x}_s)\| \notag
\\&\overset{\eqref{bd-con-err}}{\leq} N\tilde{C}  L_{fj}+\max_{ \mathbf{x}\in X}\|\phi_j( \mathbf{x})\|  \overset{\eqref{def-C}}{=}C_j.
\end{align}

By applying the optimality condition to \eqref{alg-strategy0},  using the definitions \eqref{def_breg} and \eqref{def-zeta}, we have that
\begin{equation}\label{opt-cond}
\begin{split}
&(x_i-x_{i,k+1})^T\alpha_k\big( (1+\lambda_k)(\zeta_{i,k}+F_i(x_{i,k}, N\hat{v}_{i,k+1} )\\
&-\lambda_k(\zeta_{i,k-1}+F_i(x_{i,k-1}, N\hat{v}_{i,k } )) \big)
\\& + \nabla h(x_{i,k+1}) -   \nabla h(x_{i,k})\big)\geq 0,\quad \forall x_i\in X_i.
\end{split}
\end{equation}
By setting $x_i=x_{i,k}$ in \eqref{opt-cond},   rearranging  the terms, and using the assumption that  $h $ is 1-strongly convex, we obtain
\begin{equation*}
\begin{array}{l}
\alpha_k (x_{i,k}-x_{i,k+1})^T \big( (1+\lambda_k)(\zeta_{i,k}+F_i(x_{i,k}, N\hat{v}_{i,k+1} )\\
-\lambda_k(\zeta_{i,k-1}+F_i(x_{i,k-1}, N\hat{v}_{i,k } )) \big)
\\  \geq  (x_{i,k}-x_{i,k+1})^T  \big(    \nabla h(x_{i,k})-\nabla h(x_{i,k+1})\big)
\geq    \| x_{i,k}-x_{i,k+1}\|^2.
\end{array}
\end{equation*}
Then from  \eqref{bf-F} it follows that
\begin{align*}  
\| x_{i,k}-x_{i,k+1}\| &\leq    \alpha_k   \big( (1+\lambda_k)(\| \zeta_{i,k}\|+ \|F_i(x_{i,k}, N\hat{v}_{i,k+1} \|)
+\\&\lambda_k(\| \zeta_{i,k-1}\|+\|F_i(x_{i,k-1}, N\hat{v}_{i,k } )\|) \big)  \notag \\& {\leq }   \alpha_k  \big((1+2\lambda_k)C_i  +(1+\lambda_k)\| \zeta_{i,k}\| +\lambda_k\| \zeta_{i,k-1}\| \big) .
\end{align*}
This together with \eqref{def-bdst} and \eqref{bd-consensus1}   produces \eqref{agg-bd0}.
\hfill $\Box$

\subsection{Proof of Theorem \ref{thm1}}\label{subsec-thm1}

By substituting  \eqref{bd-con-err} into the definition of
$\varepsilon_k$ in  \eqref{def-vek}, we obtain that
\begin{equation}\label{bd-varek2}
\begin{split}
\varepsilon_k& = 2\sum_{i=1}^N L_{f_i} M_i \big( \| N \hat{v}_{i,k+1}-\sigma(\mathbf{x}_k)\|
+\| N \hat{v}_{i,k }-\sigma(\mathbf{x}_{k-1})\| \big)
 \\& \leq 4N\tilde{C} \sum_{i=1}^N L_{f_i} M_i=c_e,\quad \forall k\geq 1.
\end{split} \end{equation}

With the selection of parameters, similar to the proof of \cite[Corollary 3.4]{kotsalis2020simple}, we can verify that \eqref{rela-theta},  \eqref{rela-theta2}, and \eqref{ass-theta} hold.
Furthermore,   by $c_0\geq 4$ and the simple calculations we obtain
\begin{align*}
& \theta_t (2\mu \alpha_t+1/2) >\theta_t/2=(t+c_0+1)(t+c_0 )/2,
\\& \theta_t\alpha_t^2={(t+c_0+1)(t+c_0 ) \over \mu^2(t+c_0-1)^2} \leq {2\over \mu^2},
\\& \theta_k\alpha_k\leq {2\over \mu^2\alpha_k},~ \theta_k \alpha_k \lambda_k =\theta_{k-1}\alpha_{k-1} \leq {2\over \mu^2\alpha_{k-1}}\leq {2\over \mu^2\alpha_k},
\\&\sum_{k=1}^t\theta_k \alpha_k^2 \lambda_k^2 =\sum_{k=1}^t\theta_k \alpha_k^2
{\theta_{k-1}^2\alpha_{k-1}^2 \over \theta_k^2\alpha_k^2}=
\sum_{k=1}^t  {\theta_{k-1}^2\alpha_{k-1}^2 \over \theta_k }
\\&=\sum_{k=1}^t {(k+c_0 ) (k+c_0-1 )^2 \over \mu^2  (k+c_0-2)^2(k+c_0+1) }
\leq {2t \over \mu^2}.
\end{align*}
This incorporating with \eqref{lemf5}  and \eqref{bd-varek2} implies that
\begin{equation}\label{result1}
\begin{split}
&{(t+c_0+1)(t+c_0 )\over 2}\mathbb{E} [\widehat{D}(\mathbf{x}_{t+1},\mathbf{x}^*)]
\\ & \leq   \theta_1 \mathbb{E} [\widehat{D}(\mathbf{x}_1,\mathbf{x}^*)]
+{(16t+4)\nu^2\over \mu^2} \\&+  {\theta_1  \alpha_1 \lambda_1 c_e
	+ {4\over \mu^2 } \sum_{k=1}^t {\mathbb{E} [\varepsilon_{k+1}]\over \alpha_k} }      .
\end{split} \end{equation}

In the following, we establish an upper bound of $\sum_{k=1}^t {\mathbb{E} [\varepsilon_{k+1}]\over \alpha_k} .$
Note that
\begin{equation}\label{bd-beta1}
\begin{split}
& \sum_{k=1}^{t} k \beta^{k}  =\beta  \left( \sum_{k=1}^{t} k \beta^{k-1}\right)
=\beta  \left( \sum_{k=1}^t   \beta^k\right)'=\beta\left( {\beta  - \beta^{t+1}  \over 1- \beta }\right)'
\\&={\beta\left(1-\beta^{t+1} -(t+1) \beta^t (1-\beta)\right) \over ( 1- \beta )^2 }
\leq {\beta  \over ( 1- \beta )^2 }.
\end{split}
\end{equation}
Hence  for any $t\geq 1,$
\begin{align}\label{bd-beta2}
& \sum_{k=1}^{t} {\beta^{k}  \over \alpha_k}
= \sum_{k=1}^{t} \mu(k+c_0-1)\beta^{k}=  \mu( c_0-1)\sum_{k=1}^{t}\beta^{k}+ \mu   \sum_{k=1}^{t} k\beta^{k} \notag
\\&\leq {\mu\beta (c_0-1)  \over   1- \beta  }+
{\mu\beta  \over ( 1- \beta )^2 }.
\end{align}
Also, note by $c_0\geq 4$  and \eqref{bd-beta1} that
\begin{equation}  \label{bd-beta3}\begin{split}
& \sum_{k=1}^t {1\over \alpha_{ k}} \sum_{s=1}^k \beta^{k-s}   \alpha_{ s-1}
= \sum_{k=1}^t   \sum_{s=1}^k \beta^{k-s}   { \alpha_{ s-1} \over \alpha_k}
=\sum_{k=1}^t  \sum_{s=1}^k \beta^{k-s}   { k+c_0-1\over s+c_0-2}
\\  & =\sum_{k=1}^t  \sum_{s=1}^k \beta^{k-s} \left(1+  { k-s+1\over s+c_0-2}  \right)
\leq  \sum_{k=1}^t  \sum_{s=1}^k \beta^{k-s} \left(     k-s+2  \right)
\\&   =2\sum_{k=1}^t \sum_{s=0}^{k-1} \beta^s  +\sum_{k=1}^t \sum_{s=0}^{k-1} s\beta^s
\\&  \leq  \sum_{k=1}^t  \left(  {2  \over  1- \beta   } + {\beta  \over ( 1- \beta )^2 }\right)={(2-\beta) t  \over ( 1- \beta )^2 }.
\end{split} \end{equation}
Similarly to \eqref{bd-beta3}, we have
\begin{equation}  \label{bd-beta4}\begin{split}
& \sum_{k=1}^t {1\over \alpha_{ k}} \sum_{s=1}^{k+1} \beta^{k+1-s}   \alpha_{ s-1}
=\sum_{k=1}^t  \sum_{s=1}^{k+1} \beta^{k+1-s}   \left(1+  { k-s+1\over s+c_0-2}  \right)
\\&\leq   \sum_{k=1}^t  \sum_{s=1}^{k+1} \beta^{k-s+1} \big((k-s+1)+1  \big)
\\&\leq   \sum_{k=1}^t  \left( \sum_{s=0}^{k }s \beta^{s}  +\sum_{s=0}^{k } \beta^{s}  \right)
\\&\leq  \sum_{k=1}^t  \left(  {1  \over  1- \beta   } + {\beta  \over ( 1- \beta )^2 }\right)={  t  \over   ( 1- \beta )^2 }.
\end{split} \end{equation}

Note by $c_0\geq 4$ that for any $\geq 1,$
\begin{align*}
& \lambda_k
={\theta_{k-1}\alpha_{k-1} \over \theta_k\alpha_k}={(k+c_0-1)^2\over (k+c_0+1) (k+c_0-2)}\leq 3/2.
\end{align*}
By \eqref{bd-noise}, we obtain that $\mathbb{E}[ \| \zeta_{i,k}\|]\leq \sqrt{\mathbb{E}[ \| \zeta_{i,k}\|^2]}\leq \nu_i,~ \forall k\geq 1.$ This together with  \eqref{agg-bd0}  produces
\begin{equation*}
\begin{split}
&\mathbb{E}[\|\sigma(\mathbf{x}_k) -N\hat{v}_{i,k+1}\|]
\leq \theta M_HN \beta^{k}
\\&+\theta N \sum_{s=1}^k \beta^{k-s}   \alpha_{s-1}  \sum_{i=1}^N  \big((1+2\lambda_{s-1})C_i  +(1+2\lambda_{s-1})\nu_i \big)
\\& \leq \theta M_HN \beta^{k}
+4\theta N\sum_{i=1}^N  \big( C_i  + \nu_i \big)  \sum_{s=1}^k \beta^{k-s}   \alpha_{s-1}.
\end{split}
\end{equation*}
Then by using \eqref{bd-beta2} and \eqref{bd-beta3}, we obtain that
\begin{equation}
\begin{split}
&\sum_{k=1}^t {1\over \alpha_k} \mathbb{E}[\|\sigma(\mathbf{x}_k) -N\hat{v}_{i,k+1}\|]
\\& \leq \theta M_HN \mu \beta \left({ c_0-1   \over   1- \beta  }+{1  \over ( 1- \beta )^2 }\right)
+ {4\theta N \sum_{i=1}^N   (C_i + \nu_i) (2-\beta) t  \over ( 1- \beta )^2 }.
\end{split}
\end{equation}
Similarly, by using \eqref{bd-beta2} and \eqref{bd-beta4}, we obtain that
\begin{equation}
\begin{split}
&\sum_{k=1}^t {1\over \alpha_k} \mathbb{E}[\|\sigma(\mathbf{x}_{k+1}) -N\hat{v}_{i,k+2 }\|]
\\& \leq \theta M_HN\mu \beta^2 \left({   c_0-1   \over   1- \beta  }+
{1  \over ( 1- \beta )^2 }\right)
+ {4\theta N \sum_{i=1}^N   (C_i + \nu_i)   t  \over   ( 1- \beta )^2 }  .
\end{split}
\end{equation}

By recalling  the definition of $\varepsilon_k$ in  \eqref{def-vek}, we obtain that
\begin{equation*}
\begin{split}
&\sum_{k=1}^t {\mathbb{E} [\varepsilon_{k+1}]\over \alpha_k}
=  2 \sum_{i=1}^N L_{f_i} M_i  \sum_{k=1}^t { \mathbb{E} [\| N \hat{v}_{i,k+1}-\sigma(\mathbf{x}_k)\|]\over \alpha_k}
\\&+2 \sum_{i=1}^N L_{f_i} M_i \sum_{k=1}^t { \mathbb{E} [\| N \hat{v}_{i,k+2 }-\sigma(\mathbf{x}_{k+1})\|]\over \alpha_k}
\\& \leq  2\theta M_HN \mu\beta (1+\beta)\left({  c_0-1   \over   1- \beta  }+
{ 1\over ( 1- \beta )^2 }\right)  \sum_{i=1}^N L_{f_i} M_i
\\& + t {8\theta N \sum_{i=1}^N   (C_i + \nu_i) (3-  \beta)     \over   ( 1- \beta )^2 } \sum_{i=1}^N L_{f_i} M_i=c_1+c_2t,
\end{split} \end{equation*}
where the last equality follows from  the definitions of $c_1$ and $c_2$ in \eqref{def-C1}
and \eqref{def-C2}.
This together with \eqref{result1}  and $\alpha_1\theta_1\lambda_1=\alpha_0\theta_0={ c_0 (c_0+1)\over \mu (c_0-1)}$ implies that
\begin{equation}\label{result2}
\begin{split}
&{(t+c_0+1)(t+c_0 )\over 2}\mathbb{E} [\widehat{D}(\mathbf{x}_{t+1},\mathbf{x}^*)]
\\ & \leq   \theta_1 \mathbb{E} [\widehat{D}(\mathbf{x}_1,\mathbf{x}^*)]+{(16t+4)\nu^2\over \mu^2}
\\&+  { c_0  (c_0-1)C_{e}\over \mu(c_{0}-1) }
+ {4(c_1+c_2t)\over \mu^2 }     .
\end{split} \end{equation}
Then the result follows immediately.
\hfill $\Box$

\section{Numerical Simulations}\label{sec:sim}

In this section, we validate the algorithm  performance   through   numerical simulation on the Nash-Cournot games   (see e.g., \cite{yi2019operator,lei2021distributed}).
There is a collection of $N$ factories denoted by $\mathcal{N}=\{1,\dots,N\}$ competing over $l$ markets denoted by $\mathcal{L} \triangleq\{1, \cdots, L\}$,  where  each factory $i$ needs to decide its production $x_{i,l}$ at  markets $l$. Then, the cost of production of factory $i$ is defined as $c_{i}\left(x_{i} ; \xi_{i}\right)=\left(c_{i}+\xi_{i}\right)^{T} x_{i}$, where $x_{i} \triangleq \operatorname{col}(x_{i,j})_{l=1}^{L}$, $c_{i} > 0$ is parameter for factory $i$, and $\xi_{i}$ is a random disturbance or noise with zero-mean and bounded variance.
The income of factory $i$ is $p_{l}\left(S_{l} ; \zeta_{l}\right)x_{i} $, where $S_{l}=\sum_{i=1}^{N}x_{i,l}$ denote the aggregate products of all factories delivered to market $l$. By the law of supply and demand, the price function $p_{l}$ can be represented by the reverse demand function and is defined as $p_{l}\left(S_{l} ; \zeta_{l}\right)=d_{l}+\zeta_{l}-b_{l} S_{l}$, where $d_{l}>0$, $b_{l}>0$, and  $\zeta_{l}$ is zero-mean random  disturbance or noise. Then, the factory $i$'s cost is $\psi_{i}\left(x ; \xi_{i}, \zeta_{i}\right)=c_{i}\left(x_{i} ; \xi_{i}\right)-(d+\zeta-B S)^{T}  x_{i}$, with $d=\left(d_{1}, \cdots, d_{l}\right)^{T}, B=\operatorname{diag}\left\{b_{1}, \cdots, b_{l}\right\}$, $S=\operatorname{col}\left\{S_{1},\dots,S_{l}\right\}$ and $\zeta=\operatorname{col}\left\{\zeta_{1}, \cdots, \zeta_{l}\right\} .$ Finally, factory $i$'s local optimization problem is $\min _{x_{i} \in X_{i}} f_{i}(x) \triangleq \mathbb{E}\left[\psi_{i}\left(x ; \xi_{i}, \zeta_{i}\right)\right],$ while satisfying a finite capacity constraint $X_{i}$.

It is straightforward to verify that the aforementioned  Nash-cournot example satisfy the aggregative game formulation (1) with $\sigma(x)=\sum_{i=1}^{N} x_{i}$, and $f_{i}(x, \sigma)=c_{i}^{T} x_{i}-(d-B \sigma)^{T} x_{i} .$ Then by \eqref{def-phi}, $F_{i}\left(x_{i}, z\right)=c_{i}- d+ B\left(z+ x_{i}\right)$ and $\phi_{i}(x)=$ $c_{i}- d+ B\left(\sum_{i=1}^{N}  x_{i}+ x_{i}\right)$. We can verify that  Assumptions \ref{ass-payoff}, \ref{ass-hF} and \ref{ass-noise} hold when  the random variables $\xi_{i}, \zeta_{l}, i \in \mathcal{N}, l \in \mathcal{L}$ are zero mean with bounded variance. In addition, we set $h(x)=\frac{1}{2}\|x\|_{2}^{2}$, then the Bregman divergence becomes $D_{h}(x, y)=\|x-y\|_{2}^{2}$.

Set $N=20, l=3$, and let the communication among the factories be described by an undirected time-varying graph. The graph at
each iteration is randomly drawn from a set of four graphs, whose union graph is connected.  Set the adjacency matrix $W=\left[w_{i j}\right]$, where $w_{i j}=\frac{1}{\max \left\{\left|\mathcal{N}_{i}\right|,\left|\mathcal{N}_{j}\right|\right\}}$ for any $i \neq j$ with $(i, j) \in \mathcal{E}, w_{i i}=1-\sum_{j \neq i} w_{i j}$, and $w_{i j}=0$, otherwise. Let $c_{i}$ is drawn from the uniform distribution $U[3,4]$. The  pricing parameters $d_{l}, b_{l}$ of market $l \in \mathcal{L}$ are derived from uniform distributions $U[10,10.5]$ and $U[0.5,1]$, respectively. The capacity constraint for each $x_{i,l}$ is the same as $ x_{i,l} \in [2, 10], \forall i\in \mathcal{N},l\in \mathcal{L}$.
After fixing the problem parameters $c_{i}$ and $d_{l}$,  let the random variables $\xi_{i}, i \in \mathcal{N}$ and $\zeta_{l}, l \in \mathcal{L}$ be randomly and uniformly  from $U\left[-c_{i} / 8, c_{i} / 8\right]$ and $U\left[-d_{l} / 8, d_{l} / 8\right]$, respectively.

\begin{figure}[htbp]
	\centering
	\includegraphics[width=0.36\textwidth, height=0.25\textheight]{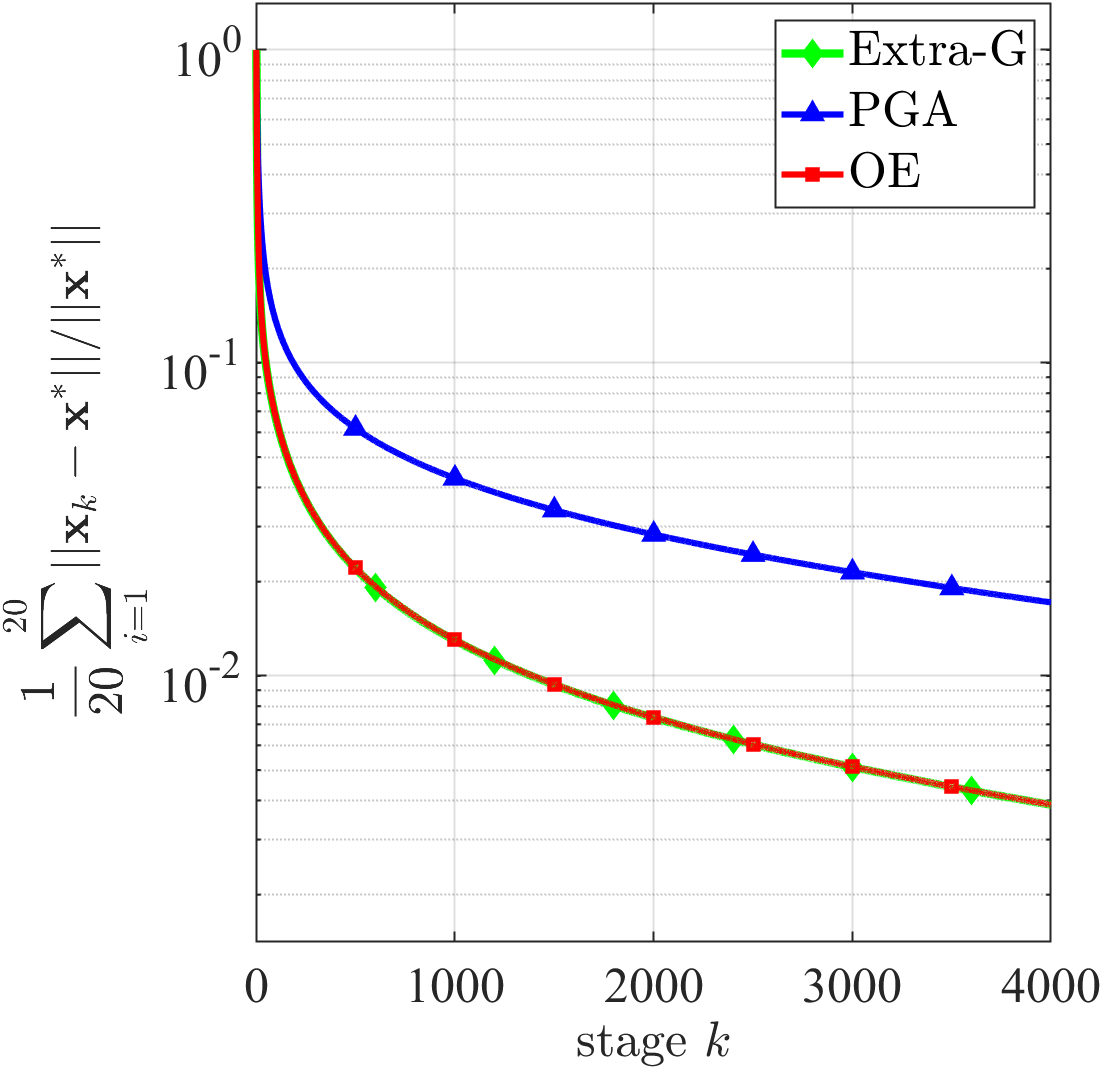}
	\caption{Comparisons of  Algorithm \ref{alg1} with the projected gradient method
 (PGA) and the  extra-gradient method  (Extra-G). The trajectories are derived by averaging with 20 sample paths.}\label{fig.com}
\end{figure}

We implement Algorithm \ref{alg1} and display the empirical results in Fig.\ref{fig.com} by averaging of 20 sampling trajectories with the same initial points, with $\alpha_k$ and $\lambda_k$ choosing as in Theorem \ref{thm1}. Besides, we compare our  algorithm with the projected gradient method \cite{koshal2016distributed} and extra-gradient method \cite{enrich2019extragradient} when applied to the stochastic aggregative game considered in this work, but the network aggregate value is still estimated with the dynamical average tracking. The projected gradient method (PGA) requires a simple projection given by
$$x_{i}^{k+1}:=\operatorname{P}_{X_{i}}\left[x_{i}^{k}-\alpha_{k} q_i(x_{i,k}, N\hat{v}_{i,k+1};\xi_{i,k})\right],$$
while the extra-gradient method (Extra-G) consists of  two steps, 
requiring two steps of projection evaluation and gradient sample at each iteration.
$$
\begin{array}{cl}
\text { (extrapolation) } & x_{i}^{k+1 / 2}=\operatorname{P}_{{X_{i}}}\left[x_{i}^{k}-\alpha_{k} q_i(x_{i,k}, N\hat{v}_{i,k+1};\xi_{i,k})\right], \\
\text { (update) } & x_{i}^{k+1}=\operatorname{P}_{{X_{i}}}\left[x_{i}^{k}-\alpha_{k} q_i(x_{i}^{k+1 / 2}, N\hat{v}_{i,k+1};\xi_{i,k})\right].
\end{array}
$$

Fig.\ref{fig.com} displays the convergence of the three algorithms, and it shows the superior convergence speed of Algorithm \ref{alg1} compared to  projected gradient method.
Though the convergence speed of Algorithm \ref{alg1} is almost the same as that of the extra-gradient method, the advantage of our method lies in that it only requires  a single projection and  one sampled gradient, greatly reducing the computing and sampling cost.

\section{Conclusions }\label{sec:con}

This paper proposes  a distributed operator extrapolation method for stochastic aggregative game based on mirror descent, and shows that the  proposed method can achieve  the optimal convergence for  the class of strongly monotone games. In addition, empirical  results demonstrate that our method indeed brings speed-ups. It is of interest to explore the algorithm convergence for monotone games,  and extend the operator extrapolation method to the other classes of network games in distributed and stochastic settings.





\bibliographystyle{IEEEtran}
\bibliography{ref}

\end{document}